# Implementation of the AC-based economic dispatch in the Russian electricity market

Mikhail Davidson
ConsultNext, SkolTech, MSU
Moscow, Russia
m.davidson@consultnext.ru

Gleb Labutin
Markets Development Department
JSC "System Operator of United Power System"
Moscow, Russia
labutin@so-ups.ru

*Abstract*—**The Russian wholesale electricity market started to operate in its present form since 2006. Its unique specificity is the AC-based economic dispatch model underlying the market applications. Solution method for the day-ahead and balancing market applications was developed and implemented. Now it is employed by the commercial operator and the system operator of the Russian energy system. It is proved to be an adequate and reliable scheduling market technology.**

## I. INTRODUCTION

In the period 2001-2008 the Russian electric power sector had undergone a substantial reform. As a result a competitive electricity market was established and a number of independent generation and distribution companies emerged. An overview of the market model is given in [3].

This required a new technology of market based short-term scheduling which was developed and implemented in 2006-2007. Since then it is exploited by the commercial market operator (Administrator of the Trading System – ATS) that runs the Day-Ahead market and the System Operator (SO UPS) that runs the Balancing market. The aim of the latter is to reflect the most recent system conditions and to trade deviations from the Day-Ahead market. Balancing market schedule is computed every one hour ahead of real time and the resulting generation station schedules have a status of dispatch instructions.

Most part of the energy system of Russia is interconnected and operates at a single frequency. The system is highly constrained and is distinctive for transferring significant amounts of power over HV backbone grid, and pattern of these transfers is changing throughout a day. A DC model widely used in many energy markets might lead to significant scheduling errors and non-transparent deviations. The AC model yields the required accuracy in modelling losses and power flows to correctly account for the constraints set up by the System Operator. For these reasons the mathematical model underlying the market applications is AC-based economic dispatch problem.

The network model used in the market applications has around 9000 nodes and 15000 buses. The time horizon is broke down into hourly intervals for which intertemporal constraints are present: ramp rate constraints and energy volume constraints (especially for hydros). Hence the resulting optimization problem is a large scale nonlinear and nonconvex program. Due to the size and nonlinearity the issues of timing and convergence reliability emerge that make developing solution algorithm a challenging problem.

The algorithm is based on sequential quadratic programming with an intensive use of parallel computations. The method is efficient for both day-ahead market and balancing market where the timing restrictions are tighter. In this paper we focus on the model setup for the balancing market for which the key difference with the day-ahead is that the load forecast is provided by the system operator while at the day-ahead consumer bids are used.

## II. MODEL SETUP AND PROBLEM FORMULATION

In this section we describe constraint types and control variables used in the model. The model includes both intertemporal and hourly constraints. The latter are referred to by the index $\tau = 1,...,T$, where T is the number of hourly intervals in the session. Also a subset of variables relating to hour $\tau$ is denoted by the corresponding index.

### A. Power balance constraints

Trading in the RF energy market concerns only active power and not reactive power. While reactive power balance equations are explicitly present in the problem formulation, optimization of reactive power is not performed in order to avoid potentially negative effects of using active power controls to enforce VQ constraints and vice-versa [9].

The values of such controls as transformer taps and generator voltage setpoints and other devices that essentially impact the reactive power balance are accepted as inputs to the problem and not varied. At the same time regulation of reactive power by generation sources is modeled through switching PV and PQ bus types depending on the maximum and minimum reactive power capacity of a source and required reactive power consumption/generation at the corresponding bus. The maximum and minimum reactive power capacity is assumed constant not depending on the active power output of a generator. In other words accounting for reactive power is similar to the load flow model.

The active power controls that are subject to optimization at the balancing market are the volume bids submitted by the generators. The load forecast is provided by the System Operator. Thus, the power balance equations connect the state variables $y_\tau$ (voltage magnitudes and angles), the active power volume bids $x_\tau$ and the fixed net bus injections $d_\tau$ representing the demand:

$$[\lambda_\tau] \quad F_\tau(y_\tau) + Bx_\tau = d_\tau, \qquad (1)$$

(here the quantity in square brackets denotes the corresponding Lagrange multipliers), $B$ is the buses - bids incidence matrix. The number of equations in (1) is twice the number of buses (excluding swing bus) reflecting active and reactive nodal power balance. We use separate notation for the swing bus balance:

$$[\lambda_0] \quad F_{0\tau}(y_\tau) = 0. \qquad (2)$$

Since we don't optimize for reactive power controls the reactive power swing bus equation is slack and (2) represents only active power balance (the right hand side is zero to simplify the notation).

### B. Power flow constraints

Power flow constraints represent a set of most likely contingency scenarios under given operating conditions modelled and precomputed by the System Operator. For every such constraint SO defines a power flow limit that guarantees system security after any single contingency, and these are the inputs to the model:

$$[\sigma_\tau] \quad G_\tau(y_\tau) \le \bar{f}_\tau. \qquad (3)$$

### C. Intertemporal constraints

These include ramp rate and fuel generator constraints. The latter model the daily water supply of hydros to be optimally allocated among hourly intervals. Both ramp rate and fuel constraints are linear and link only volume variables. We denote these constraints by a single inequality (including upper and lower bounds):

$$ITx \le b_{IT} \qquad (4)$$

where $x = (x_1, ..., x_T)^T$.

### D. Bounds for variables

For generation volumes $x_\tau$ there are upper and lower bounds submitted by the generators in their bids

$$lb_\tau \le x_\tau \le ub_\tau. \qquad (5)$$

We don't apply bounds on voltage magnitudes since as said above the controls affecting reactive power and voltage magnitudes are not optimized.

### E. Problem formulation

The objective function used at the RF balancing market is minimum of the cost of energy produced by the generators according to their price bids. Thus, the problem formulation is

$$\sum_\tau c_\tau^T x_\tau \to \min \qquad (6)$$

subject to constraints (1) – (5). This is a nonlinear, nonconvex and large-scale optimization problem. The size of the problem is of order 5e+5 variables and constraints. Modeling PV-PQ bus switching brings a discrete discontinuous nature to the problem.

### III. THE ALGORITHM

SQP methods proved to be efficient in application to general nonlinear problems [5]. In particular, successful application of SNOPT [4] to power flow optimization is reported in [2]. Here we describe a variant of SQP algorithm used in the Russian energy market.

For a given pair $(y, x)$ and the corresponding vector of Lagrange multipliers $(\lambda, \sigma)$ generic form SQP sub-problem as applied to (5) appears as follows

$$\sum_\tau \left(c_\tau^T \Delta x_\tau + 1/2 \Delta y_\tau^T H_\tau(\lambda_\tau, \sigma_\tau y_\tau) \Delta y_\tau \right) \to \min \qquad (7)$$

subject to

$$DF_\tau(y_\tau)\Delta y_\tau + B_\tau \Delta x_\tau = r_\tau^F \qquad (8)$$

$$DF_{0\tau}(y_\tau)\Delta y_\tau = r_\tau^{0\tau} \qquad (9)$$

$$DG_\tau(y_\tau)\Delta y_\tau \le r_\tau^G \qquad (10)$$

$$IT\Delta x \le r^{IT} \qquad (11)$$

$$lb_\tau \le x_\tau + \Delta x_\tau \le ub_\tau \qquad (12)$$

where (8) – (10) are linearization of constraints (1) – (3) around $y$, and

$$H_\tau(\lambda_\tau, \sigma_\tau y_\tau) = \sum_i \left(\lambda_\tau^i \nabla_{yy}^2 F_{\tau i}(y_\tau)\right) + \sum_j \left(\sigma_\tau^j \nabla_{yy}^2 G_{\tau j}(y_\tau)\right)$$

is the element of the Hessian of the Lagrange function of (6) corresponding to hour $\tau$. Only Hessians of constraints (1) – (3) enter the formula for $H_\tau$ (let the first sum contain Hessian terms of both equations (1) and (2)).

Despite the size of the problem is large, the number of degrees of freedom is much smaller. As long as the state variables $y$ stay away from the boundary of the solvability region the Jacobian $DF_\tau(y_\tau)$ is nonsingular and the states $y$ could be seen as implicit functions of controls $x$, and the state increments could be expressed via the control increments:

$$\Delta y_\tau = DF_\tau^{-1}(y_\tau)r_\tau^F - DF_\tau^{-1}(y_\tau)B_\tau \Delta x_\tau. \qquad (13)$$

Clearly dimension of the controls (the number of generator bids) is much smaller than the number of buses. Note that

$DF_\tau(y_\tau)$ is sparse and sparse LU factorization could be used to efficiently perform the required computations. The reduced gradients of $F_{0\tau}$ and $G_\tau$ could be readily computed. Despite the reduced gradients does not keep sparse structure, the size of the corresponding matrices is generally small. Dimension of $F_{0\tau}$ is the number of swing buses (islands) in the system. Most part of the system operates as a single whole, however temporarily a part may operate separately due to grid maintenance. In practice the number of islands does not exceed 3. Flow constraints (3) are managed using an extended active set strategy: inequality (10) contains linearized constraints (3) binding at the current iterate $y_\tau$ or at any of the previous iterates. At a cost of slightly expanding the active set we avoid oscillations and stabilize the convergence. In our experience the number of the elements in the extended active set rarely exceeds 20 per hourly interval. Hence computation of the reduced gradients for (2) and (3) is cheap.

Additional computational savings could be achieved by employing the fact that the structure of the Jacobian does not change from iteration to iteration: the symbolic analysis of the matrix structure required for the sparse LU can be performed only once. However, this holds only until there is PV-PQ bus type switching.

Computation of the reduced Hessian is the most time consuming. The reduced Hessian has the form

$$\tilde{H}_\tau(\lambda_\tau, \sigma_\tau y_\tau) = B_\tau^T (DF_\tau^{-1})^T H_\tau DF_\tau^{-1} B_\tau.$$

The L, U factors of $DF_\tau$ are not sparse and so is $\tilde{H}_\tau$. While dimension of $\tilde{H}_\tau$ is much smaller than that of $H_\tau$, we in addition take into account that the columns of $B_\tau$ are identical for generators located at the same bus. For computing the reduced Hessian it is sufficient to keep only one column per bus which reduces the time for computing $\tilde{H}_\tau$ up to three times.

In order that the solution to the SQP sub-problem be well defined $\tilde{H}_\tau$ has to be positive semidefinite. In general this does not hold. In the algorithm we use a positive semidefinite approximation to the reduced Hessian denoted by $\tilde{H}_\tau^+$.

In reduced form the SQP sub-problem (7) will appear as follows

$$\sum_\tau \left( \tilde{c}_\tau^T \Delta x_\tau + 1/2 \Delta x_\tau^T \tilde{H}_\tau^+ \Delta x_\tau \right) \to \min \qquad (14)$$

subject to constraints

$$L_{0\tau} \Delta x_\tau = r_\tau^L \qquad (15)$$
$$S_\tau \Delta x_\tau \le r_\tau^S \qquad (16)$$

and (11), (12) that remain unchanged. Here $L_{0\tau}$ and $S_\tau$ are the reduced gradients of $F_{0\tau}$ and $G_\tau$, respectively, $\tilde{c}_\tau$ is the combination of $c_\tau$ and terms that appear after substitution (13):

$$\tilde{c}_\tau = c_\tau - (r_\tau^F)^T (DF_\tau^{-1})^T H_\tau DF_\tau^{-1} \qquad (17)$$

Note that $L_{0\tau}$ is loss coefficient vector and $S_\tau$ is power transfer distribution factor matrix representing active loss and power flow constraint sensitivities to bus injections. In particular

$$L_{0\tau} = 1 + \partial L_\tau / \partial x_\tau^a \qquad (18)$$

where on the right hand side there is the gradient of the active power loss function with respect to active power nodal injections.

The optimal solution $\Delta x_\tau^+$ of (14) is used to update the current iterate $x_\tau$. State variables are updated via (13). To update the dual multipliers $\lambda_\tau$ note that from the first order necessary conditions for (7) one has:

$$\lambda_\tau^T + \lambda_{0\tau} DF_{0\tau} DF_\tau^{-1} + \sigma_\tau^T DG_\tau DF_\tau^{-1} = 0 \qquad (19)$$

(for the dual multipliers in (7) we use the same notation as for (6)). This yields the well known decomposition of the locational marginal prices (LMPs) through marginal cost of losses and marginal cost of transmission constraints [1].

Thus (14) is solved for $\Delta x_\tau^+$, $\lambda_{0\tau}^+, \sigma_\tau^+$ that are used to update the current iterate and to build a new quadratic sub-problem. The following specifications complete the definition of the SQP algorithm:

- Initial approximation is always chosen to be all zero angles and unit voltages ("flat start"). The multipliers $\lambda_\tau$ are initialized according to (19) with $\lambda_{0\tau}$ being some positive proxy of the system price. The multipliers $\sigma_\tau$ are zero as transmission constraints are not tight at flat start state.

- Since the algorithms starts from an infeasible state, at the first stage the step size is chosen in order to minimize the error of the power balance nodal equations; if at some iteration the feasibility error is within the specified tolerance the step size is chosen to minimize a combination of the objective and the feasibility error.

- The stopping criterion is defined using the error of the first order necessary conditions for (6).

This breaks down into checking that the current state is feasible to a given tolerance and that the norm of $\tilde{H}_\tau^+ \Delta x_\tau$ is

within the optimality tolerance. The first order optimality conditions for the SQP sub-problem (14) imply

$$\tilde{c}_\tau + \tilde{H}_\tau^+ \Delta x_\tau + L_{0\tau}^T \lambda_{0\tau} + S_{0\tau}^T \sigma_\tau + \pi_\tau = 0 \qquad (20)$$

where $\pi_\tau$ denote the aggregate dual term related to constraints (11) and (12). For every unit $\pi_\tau$ is a marginal profit/loss of the unit given the system price and the unit price bid. Since the state is feasible up to a given tolerance, the residual $r_\tau^F$ is small and then from (17) $\tilde{c}_\tau$ is close to $c_\tau$, and (20) is close to

$$c_\tau + \pi_\tau = L_{0\tau}^T \lambda_{0\tau} + S_{0\tau}^T \sigma_\tau \qquad (21)$$

The complementarity conditions in (14) provide that unit's marginal profit $\pi_\tau$ is greater or equal zero as long as the unit operates at its maximum capacity and less or equal zero if it is at its minimum. Then (21) yields the equilibrium condition for every unit at given system price at the right hand side of (21).

Since energy pricing in the Russian energy market is based on LMPs, conditions (20) make the resulting schedule and prices clear, transparent and justifiable for the market participants. Conditions (20) have dimension of national currency units per MWh. Hence the choice of the optimality tolerance could be linked to the minimum fraction of the currency unit (for example, one hundredth) up to which the price is stated for accounting purposes.

Finally in this section we note that the computations required to evaluate the constraints and the reduced gradients and the hessian matrix are independent for each time interval which allows using parallelization. This leads to substantial time savings as parallelization covers the most time consuming computations.

The reduced problem has the dimensions an order of magnitude less than the original problem and can be efficiently solved by standard QP-solvers.

## IV. CONVERGENCE

Conditions that guarantee local convergence of the SQP algorithm are studied in a number of papers, see [5] for a review. In [7] the local quadratic rate of convergence to a solution $x^*$ is proved under the following assumptions:

*a) Second order necessary optimality conditions hold at the solution.*

*b) The gradients of the active constraints at the solution are linearly independent.*

*c) Strict complementarity slackness holds.*

As shown in [10] last two conditions could be relaxed. Instead of b) it is sufficient that Mangasarian-Fromovitz [6] constraint qualification holds, i.e.

*d) the gradients of equalities (in our case (15)) are linearly independent, and there is a point where inequalities (11), (12), and (16) are strict.*

While the original condition of linear independence of the gradients in fact implies the uniqueness of optimal dual solution, in the relaxed form it is sufficient to require just an existence of such dual solution for which the strict complementarity holds.

In order to analyze the structural properties of optimization problem arising in the Russian energy market from the point of view of the above sufficient conditions we tried to select a representative set of sessions for an empirical study. The set covers the first half of 2017 and includes two sessions for every month that correspond to a working day and a weekend. The total number of selected sessions is thus twelve each containing 24 hourly intervals. For each session we took a solution satisfying the stopping criterion to which the method converged and analyzed the above conditions.

It turned out that b) is not satisfied for any of the sessions selected while d) holds for all of them. We were not able to verify the relaxed form of c), but for the dual solution obtained as a result of the iteration process the strict complementarity does not hold. However if one takes not all the active constraints but only nonlinear active constraints (whose linearized form is (15) and (16)) then c) is satisfied for them in all cases considered, and moreover, the gradients of these constraints are linear independent. In general this fact by itself does not imply uniqueness of the corresponding dual multipliers. But for the problem in question these multipliers define the nodal prices. Throughout the test set considered the nodal prices were unique for the primal solution obtained. This does not imply, however, the uniqueness of the whole vector of the dual solution. The linear constraints (11) and (12) are box constraints or ramp-rate constraints for units, or daily energy volume constraints for stations. Typical examples of degeneracy here might be bids with the same price at the same bus or units with ramp rate exactly equal to the MWt range between minimum and maximum capacity of the unit. Some of such degeneracies could be preprocessed but in our experience they don't essentially affect the convergence since they don't affect the nodal prices. Most important is the uniqueness of the nodal prices as they enter the expression for the Hessian, and are critical for the convergence.

The Hessian $H_\tau(\lambda_\tau, \sigma_\tau y_\tau)$ is generally not positive definite and condition a) does not hold. However, it is positive definite at "flat start" state. Let $\lambda_\tau^a$ denote dual multipliers to active power balance constraints, and $\lambda_\tau^r$ denote dual multipliers to reactive power balance constraints (for a particular $\tau$). Note that power flow constraints prices $\sigma_\tau$ are zero at "flat start". Then, the multipliers $\lambda_\tau^a$, which are the nodal prices, can be expressed as

$$\lambda_\tau^a = \lambda_{0\tau}(1 + \partial L_\tau / \partial x_\tau^a),$$

(see (18), (19)). Similarly

$$\lambda_\tau^r = \lambda_{0\tau} \partial L_\tau / \partial x_\tau^r,$$

where $\partial L_\tau / \partial x_\tau^r$ is the gradient of the active power loss function with respect to reactive power nodal injections.

At "flat start" the components of the gradient of the active power loss function are zero, and hence nodal prices $\lambda_\tau^a$ are equal for all nodes and vector $\lambda_\tau^r$ is zero. In such a case $H_\tau$ coincides (up to a multiplier) with the Hessian of the active power loss function (as a functions of state variables) which is positive definite if line reactance in the system is positive. In some neighborhood of the "flat start" state corresponding to "normal" operating conditions (where the components of the loss function gradient are sufficiently small) the hessian $H_\tau$ is close to convex and the convergence is fast. In abnormal cases when operating point approaches boundary the loss function gradient tends to infinity and convergence is violated.

Summarizing the results of the analysis we conclude that while from a formal point of view the sufficient conditions are not fully satisfied, the violations are not significant for the convergence unless the operating point is close to the boundary of the solvability set. The latter usually shows up in some local fragment of the system which is the source of the problem, and it is critically important to develop mechanisms that help to identify such "weak" locations.

## V. CONCLUSIONS

For the period of operation, the model proved to be an adequate short term scheduling market technology. The locational marginal pricing mechanism provides economic signals to participants that are concordant with the dispatch. The accuracy of the model could be illustrated by the amount of dispatch instructions issued beyond the computed schedule. According to the official report of the SO UPS for the year 2016 [8] this is within 1% of the scheduled amounts. Also according to the ATS data (www.rosenergo.ru), the corresponding financial imbalance (or "make-whole" payment) constitutes around 1% of the treaded volume, hence the resulting wholesale price is only slightly blurred by the regulatory component.

The technology has also proved to be reliable. However since the scheduling is time critical, a possibility for the dispatch personnel to analyze the resulting electrical regimes (especially in the case of problems or convergence breakdowns) is limited. Further research and work are aimed at developing algorithms for identification of weak locations in the energy system that may be a cause of divergence and providing the corresponding information to operation personnel.